\documentclass{article}

\usepackage{amsfonts}
\usepackage{amsthm}
\usepackage{amsmath}
\usepackage{amssymb}
\usepackage{mathtools} 
\usepackage{authblk} 
\usepackage{xcolor}

\usepackage[top=1in, bottom=1in, left=1in, right=1in]{geometry}

\newtheorem{thm}{Theorem}[section]
\newtheorem{lem}{Lemma}[section]

\theoremstyle{definition}

\begin{document}


\title{Solving PDEs of fractional order using the unified transform method}

\date{}


\author[1]{Arran Fernandez\thanks{Corresponding author. Email: \texttt{af454@cam.ac.uk}}}
\author[2,3]{Dumitru Baleanu\thanks{Email: \texttt{dumitru@cankaya.edu.tr}}}
\author[1,4]{Athanassios S. Fokas\thanks{Email: \texttt{T.Fokas@damtp.cam.ac.uk}}}

\affil[1]{{\small Department of Applied Mathematics and Theoretical Physics, University of Cambridge, Wilberforce Road, CB3 0WA, United Kingdom}}
\affil[2]{{\small Department of Mathematics, Cankaya University, 06530 Balgat, Ankara, Turkey}}
\affil[3]{{\small Institute of Space Sciences, Magurele-Bucharest, Romania}}
\affil[4]{{\small Viterbi School of Engineering, University of Southern California, Los Angeles, CA 90089, USA}}



\maketitle


\begin{abstract} We consider the unified transform method, also known as the Fokas method, for solving partial differential equations. We adapt and modify the methodology, incorporating new ideas where necessary, in order to apply it to solve a large class of partial differential equations of fractional order. We demonstrate the applicability of the method by implementing it to solve a model fractional problem.
\end{abstract}











\section{Introduction}

\subsection{Fractional calculus} \label{intro:FC}

Fractional calculus is a rapidly growing branch of science which is concerned with differentiation and integration to orders beyond the integers. This classical field of study, which has been considered by several great mathematicians including Leibniz, Riemann, Weyl, and Hardy \cite{miller,samko}, has become popular in recent decades due to its many applications in various fields of science \cite{hilfer,baleanu}.

There are multiple ways in which fractional derivatives and integrals can be defined; these are not all equivalent to each other, but each of them has its own advantages and disadvantages \cite{samko,atangana}. Here we shall use the classical Riemann--Liouville model of fractional calculus, in which fractional integrals are defined by
\begin{equation}
\label{RLdef:int}
\prescript{}{a}I_x^{\alpha}q(x)=\frac{1}{\Gamma(\alpha)}\int_a^x(x-\xi)^{\alpha-1}q(\xi)\,\mathrm{d}\xi,\quad\quad\mathrm{Re}(\alpha)>0,
\end{equation}
and fractional derivatives are defined by
\begin{equation}
\label{RLdef:deriv}
\prescript{}{a}D_x^{\alpha}q(x)=\frac{\mathrm{d}^m}{\mathrm{d}x^m}\left(\prescript{}{a}I_x^{m-\alpha}q(x)\right),m=\lfloor\mathrm{Re}(\alpha)\rfloor+1,\quad\quad\mathrm{Re}(\alpha)\geq0.
\end{equation}
In general, fractional derivatives and integrals can be seen as two cases of a single class of operators called \textit{differintegrals}, where the only difference between differentiation and integration is in the sign of the real part of the order variable. We may use any of the following equivalent notations for the $\alpha$th fractional differintegral of a function $q$: \[\prescript{}{a}D_x^{\alpha}q(x),\quad D_{a+}^{\alpha}q(x),\quad \prescript{}{a}I_x^{-\alpha}q(x),\quad \frac{\mathrm{d}^{\alpha}q}{\mathrm{d}(x-a)^{\alpha}}.\]

The parameter $a$ in each case is a constant of differintegration, which is usually taken to be either $0$ or $-\infty$. This is analogous to a constant of integration in classical calculus, namely, a lower bound which must be specified \textit{a priori} for an unambiguous definition. Different values of $a$ may be required in different contexts in order to obtain appropriate results; for example, we have the following formulae for fractional derivatives of power functions and exponential functions:
\begin{alignat}{2}
\label{RL:power}\prescript{}{0}D^{\alpha}_x(x^n)&=\frac{\Gamma(n+1)}{\Gamma(n-\alpha+1)}x^{n-\alpha}, &&\quad\quad\mathrm{Re}(n)>-1; \\
\label{RL:exp}\prescript{}{-\infty}D_x^{\alpha}(e^{nx})&=n^{\alpha}e^{nx},&&\quad\quad n\not\in(-\infty,0].
\end{alignat}
In a similar manner, we can define fractional differintegrals with an upper limit of differintegration instead of a lower limit:
\begin{alignat}{2}
\label{RLdef:int2} \prescript{}{x}I_b^{\alpha}q(x)&=\frac{1}{\Gamma(\alpha)}\int_x^b(\xi-x)^{\alpha-1}q(\xi)\,\mathrm{d}\xi,&&\quad\quad\mathrm{Re}(\alpha)>0; \\
\label{RLdef:deriv2} \prescript{}{x}D_b^{\alpha}q(x)&=(-1)^m\frac{\mathrm{d}^m}{\mathrm{d}x^m}\left(\prescript{}{x}I_b^{m-\alpha}q(x)\right),m=\lfloor\mathrm{Re}(\alpha)\rfloor+1,&&\quad\quad\mathrm{Re}(\alpha)\geq0.
\end{alignat}
For fractional differintegrals of this type, we may use any of the following equivalent notations: \[\prescript{}{x}D_b^{\alpha}q(x),\quad D_{b-}^{\alpha}q(x),\quad \prescript{}{x}I_b^{-\alpha}q(x).\]

It is important to note that fractional derivatives, when composed with each other, do not always yield a further fractional derivative of the appropriate order. We have the following results governing such behaviour \cite{miller,samko}.

\begin{lem}[Composition rules]
\label{ComposeDerivatives}
Fractional differintegrals of fractional integrals behave as expected, namely \[\prescript{}{a}I_x^{\alpha}\left(\prescript{}{a}I_x^{\beta}q(x)\right)=\prescript{}{a}I_x^{\alpha+\beta}q(x),\quad\quad\alpha,\beta\in\mathbb{C},\,\mathrm{Re}(\beta)>0,\] provided the relevant differintegrals exist. However, fractional differintegrals of derivatives may not behave as expected: we have \[
\prescript{}{a}D_x^{\alpha}\big(\prescript{}{a}D_x^nq(x)\big)=\prescript{}{a}D_x^{\alpha+n}q(x)-\sum_{k=1}^n\frac{(x-a)^{-\alpha-k}}{\Gamma(-\alpha-k+1)}q^{(n-k)}(a),\quad\quad n\in\mathbb{N},\,\alpha\in\mathbb{C},\] provided the relevant differintegrals exist.
\end{lem}

For this reason, it is possible to obtain a different, non-equivalent, definition of fractional derivatives by exchanging the operations of fractional integration and standard differentiation in the definition \eqref{RLdef:deriv} of Riemann--Liouville fractional derivatives. This is called the Caputo definition, and we denote it as follows:
\begin{equation}
\label{CAPdef}
\prescript{C}{a}D_x^{\alpha}q(x)=\prescript{}{a}I_x^{m-\alpha}\left(\frac{\mathrm{d}^m}{\mathrm{d}x^m}q(x)\right),m=\lfloor\mathrm{Re}(\alpha)\rfloor+1,\quad\quad\mathrm{Re}(\alpha)\geq0.
\end{equation}

We also state the following fractional generalisation of the integration by parts law, whose proof can be found in \cite{agrawal}, and which we shall need to use in \S\ref{S:GRfour} below.

\begin{lem}[Integration by parts]
\label{IbPlemma}
Let $[a,b]$ be an interval in $\mathbb{R}$ and $\alpha$ be a complex number with $\mathrm{Re}(\alpha)>0$. We have
\begin{equation}
\label{IbP:fractional}
\int_a^bf(x)\cdot \prescript{}{a}D_{x}^{\alpha}g(x)\,\mathrm{d}x=\int_a^bg(x)\cdot \prescript{C}{x}D_{b}^{\alpha}f(x)\,\mathrm{d}x-\sum_{j=0}^{n-1}\Big[(-1)^{n+j}\prescript{}{a}D_{x}^{\alpha-n+j}g(x)\cdot \prescript{}{a}D_{x}^{n-j-1}f(x)\Big]_a^b,
\end{equation}
provided the relevant differintegrals exist, where $n:=\lceil\mathrm{Re}(\alpha)\rceil$.
\end{lem}

One particularly important subfield of fractional calculus is the study of fractional differential equations. Entire textbooks have been written specifically on this topic \cite{kilbas,podlubny}, but the state of the field is still much less advanced than that of integer-order differential equations. Even ordinary differential equations, when taken to fractional orders, can still be challenging to solve \cite{area,li}, and fractional partial differential equations are of course even harder. Many of the new advances in this field are made by adapting methods from classical calculus so that they can be applied to fractional equations as well: see \cite{baleanu,bin,fernandez,podlubny2,yang} for examples of such methods. In this paper, we examine the extension to fractional partial differential equations of another powerful method for solving classical partial differential equations: namely, the unified transform method.

\subsection{The unified transform method} \label{intro:UTM}

The \textbf{unified transform method}, also known as the \textbf{Fokas method}, for solving partial differential equations is a novel technique due to the third author \cite{fokas1}. It involves integral transforms with respect to both spatial and temporal variables, where both types of transform are applied simultaneously. It is more widely applicable than classical transform methods, and can be used in several contexts where the classical transforms fail, including certain classes of evolution PDEs formulated on the half-line and the finite interval. Most importantly, the method is \textit{constructive}, generating explicit solutions in integral form. See \cite{fokas2} and \cite{fokas3} for more detail about this method and its applicability.

One important context where this method can be applied is in solving equations of the form
\begin{equation}
\label{ex:PDE}
q_t+w\big(\hspace{-0.1cm}-i\tfrac{\partial}{\partial x}\big)q=0,\quad\quad x\in(0,\infty),t\in(0,T),
\end{equation}
where $w$ is a polynomial function such that $\mathrm{Re}(w(k))\geqslant0\;\forall k\in\mathbb{R}$, with initial condition $q(x,0)=q_0(x)$ (for some known function $q_0$) and appropriate boundary conditions to be fixed later. For equations of the above form, the unified transform method works as follows.

\begin{description}
\setlength{\itemindent}{0cm}

\item[Divergence form] Introducing an exponential term enables us to write the given PDE as a family of PDEs in divergence form, parametrised by a new complex variable $k$. Specifically, the PDE can be rewritten in the form
\begin{equation}
\label{ex:divergence}
\big(e^{-ikx+w(k)t}q\big)_t=\big(e^{-ikx+w(k)t}Q\big)_x,
\end{equation}
\noindent where the function $Q(x,t,k)$ must satisfy \[\left(\frac{\partial}{\partial x}-ik\right)Q=\left(w(k)-w\left(-i\frac{\partial}{\partial x}\right)\right)q.\] Since $w$ is a polynomial, $Q$ can be defined as a finite series:
\begin{equation}
\label{ex:Qdef}
Q(x,t,k)=i\left(\frac{w(k)-w(l)}{k-l}\right)\bigg|_{l=-i\frac{\partial}{\partial x}}(q)=\sum\limits_{j=0}^{n-1}c_j(k)\frac{\partial^jq}{\partial x^j},
\end{equation}
\noindent for some complex polynomials $c_0,c_1,\dots,c_{n-1}$.

\item[Global relation] Re-expressed in divergence form, the PDE can now be integrated with respect to both $x$ and $t$. Specifically, we first substitute $\tau$ for $t$ and then apply the two operators $\int_0^{\infty}\mathrm{d}x$ and $\int_0^t\mathrm{d}\tau$ to the divergence form \eqref{ex:divergence}. On the left-hand side, the $t$-derivative disappears and we get an $x$-integral transform, which turns out to be the Fourier transform. On the right-hand side, the $x$-derivative disappears and we get a $t$-integral transform, which is a relative of the Fourier transform but considerably more complicated. The resulting identity is called the \textbf{global relation} and is central to the applicability of the unified transform method:
\begin{equation}
\label{ex:GR}
e^{w(k)t}\hat{q}(k,t)=\hat{q}_0(k)-\tilde{g}(k,t),\quad\quad\mathrm{Im}(k)<0,t\in(0,T).
\end{equation}
Here $\hat{}\;$ denotes the Fourier transform and $\tilde{g}(k,t)$ is a more complicated function, involving $t$-integral transforms with kernel $e^{w(k)t}$ applied to the functions $\frac{\partial^jq}{\partial x^j}(0,t)$ for $j=0,1,\dots,n-1$.

\item[Integral formula] The global relation \eqref{ex:GR} is an expression for the $x$-Fourier transform of $q(x,t)$ in terms of various initial and boundary values. Applying the Fourier inversion theorem yields an integral expression for $q(x,t)$ itself:
\begin{equation}
\label{ex:qint}
q(x,t)=\frac{1}{2\pi}\int_{-\infty}^{\infty}e^{ikx-w(k)t}\hat{q}_0(k)\,\mathrm{d}k-\frac{1}{2\pi}\int_{-\infty}^{\infty}e^{ikx-w(k)t}\tilde{g}(k,t)\,\mathrm{d}k.
\end{equation}
Now by Cauchy's theorem, the contour used for the second integral can be deformed from the real line to the boundary of the domain \[D^+\coloneqq\{k\in\mathbb{C} : \mathrm{Re}(w(k))<0,\mathrm{Im}(k)>0\}.\] The choice of the domain $D^+$ is motivated by considerations of exponential growth and decay: both of the exponential terms $e^{ikx}$ and $e^{-w(k)t}$ appearing in the integrand should decay as $k\rightarrow\infty$ in the regions through which the contour is deformed.

\item[The final result] The formula \eqref{ex:qint} is not the final form of the solution: it expresses $q$ in terms of the given initial condition $q_0$ and boundary values consisting of the functions $\frac{\partial^jq}{\partial x^j}(0,t)$ for $0\leqslant j<n$, which are considerably more boundary values than the number of boundary conditions necessary for the problem to be well-posed. The final step of the unified transform method involves substituting the global relation \eqref{ex:GR} into the equation \eqref{ex:qint} in order to eliminate the unknown boundary values.

The global relation holds for $\mathrm{Im}(k)<0$, while the contour of integration $\partial D^+$ is contained in the upper half plane, so some substitutions will have to be made. We replace $k$ in \eqref{ex:GR} by $\nu(k)$, where $\nu$ is a $w$-preserving function ($w(\nu(k))\equiv w(k)$) mapping $\partial D^+$ into the lower half plane, and then we use the resulting identity in \eqref{ex:qint}.

In general, there exist several possible functions $\nu$, and the identities resulting from them give several simultaneous equations in a similar form to \eqref{ex:qint}. From these equations, the unknown boundary values can then be eliminated.
\end{description}

For example, let us examine how the method would be applied to the following third-order PDE on the half-line:

\begin{align*}
q_t+q_{xxx}=0,\quad\quad &x\in(0,\infty), t\in(0,T);\\
q(x,0)=q_0(x),\quad\quad &x\in(0,\infty); \\
q(0,t)=g_0(t),\quad\quad &t\in(0,T).
\end{align*}

Let us introduce the notations $g_1(t)=q_x(0,t)$ and $g_2(t)=q_{xx}(0,t)$, where both $g_1$ and $g_2$ must be eliminated from the integral representation at a later stage.

\begin{description}
\setlength{\itemindent}{0cm}

\item[Divergence form] Here $w(k)=-ik^3$, so \eqref{ex:Qdef} yields \[Q=\Big(\frac{k^3-l^3}{k-l}\Big)\Big|_{l=-i\frac{\partial}{\partial x}}(q)=-q_{xx}-ikq_x+k^2q.\]

\item[Global relation] \[e^{-ik^3t}\hat{q}(k,t)=\hat{q}_0(k)-k^2\tilde{g}_0(-ik^3,t)+ik\tilde{g}_1(-ik^3,t)+\tilde{g}_2(-ik^3,t).\]

\item[Integral formula]
\begin{multline*}
q(x,t)=\frac{1}{2\pi}\int_{-\infty}^{\infty}e^{ikx+ik^3t}\hat{q}_0(k)\,\mathrm{d}k \\ -\frac{1}{2\pi}\int_{\partial D^+}e^{ikx+ik^3t}\big(k^2\tilde{g}_0(-ik^3,t)-ik\tilde{g}_1(-ik^3,t)-\tilde{g}_2(-ik^3,t)\big)\,\mathrm{d}k,
\end{multline*}
where the domain $D^+$ in this case is the infinite sector $\{k\in\mathbb{C}:\tfrac{\pi}{3}<\arg(k)<\tfrac{2\pi}{3}\}$.

\item[The final result] We need transformations $\nu$ such that $\nu(k)^3=k^3$. Using $\nu(k)=\omega k$ and $\nu(k)=\omega^2k$, where $\omega$ is the cube root of unity, transforms the global relation to the following equations valid for $k\in\partial D^+$:
\begin{align*}
e^{-ik^3t}\hat{q}(\omega k,t)&=\hat{q}_0(\omega k)-\omega^2k^2\tilde{g}_0+i\omega k\tilde{g}_1+\tilde{g}_2; \\
e^{-ik^3t}\hat{q}(\omega^2k,t)&=\hat{q}_0(\omega^2k)-\omega k^2\tilde{g}_0+i\omega^2k\tilde{g}_1+\tilde{g}_2(-ik^3,t).
\end{align*}
Some simple algebraic manipulations yield the final result, namely a formula for $q$ which depends only on $q_0$ and $g_0$ but not on $g_1$ and $g_2$: \[q(x,t)=\frac{1}{2\pi}\int_{-\infty}^{\infty}e^{ikx+ik^3t}\hat{q}_0(k)\,\mathrm{d}k-\frac{1}{2\pi}\int_{\partial D^+}e^{ikx+ik^3t}\big(3k^2\tilde{g}_0(-ik^3,t)-\omega\hat{q}_0(\omega k)-\omega^2\hat{q}_0(\omega^2k)\big)\,\mathrm{d}k.\]
\end{description}

A major advantage of the unified transform method is that, unlike traditional transform methods, it constructs representations of solutions which are always uniformly convergent at the boundaries of the domain. This makes it straightforward to verify that the solution does indeed satisfy the appropriate boundary conditions \cite{millersmith}. The unified transform method is also well-suited to solving a wide variety of different boundary value problems: Dirichlet, Neumann, Robin, and even non-separable boundary value problems. Furthermore, the unified transform method gives rise to simple numerical techniques for computing solutions; the numerical aspects have been explored extensively in the literature, see for example \cite{ashton,crooks,fornberg}.

\section{Applying the method to linear fractional PDEs} \label{S+:main}

There exist some earlier works analysing fractional PDEs on the half-line using methods similar to the unified transform method, e.g. the recent work of Arciga et al \cite{arciga2,arciga3,arciga4} and some papers of Kaikina \cite{kaikina1,kaikina2}. However, some of these works have used other fractional models than the classical Riemann--Liouville one -- such as the Riesz, Caputo, or Abel fractional derivatives \cite{samko} -- or considered a narrower class of PDEs than that which we shall analyse here, while others have used more complicated methods than the unified transform method. There are also issues concerning branch cuts, which naturally arise when non-integer power functions are introduced; some of the above-cited papers have skirted around these issues, while here we address them carefully and consider what deformations of contours in the complex plane are permissible when branch cuts are excluded from the domain.

Our work is both rigorous and elementary, discussing and avoiding several potential pitfalls and constructing a clear explicit algorithm for solving a large class of linear fractional PDEs on the half-line. In what follows we will address the problems which arise when a simple polynomial is replaced by a general linear combination of power functions, and we will also consider and resolve the issue of branch cuts arising from the complex power functions involved in the analysis.

\subsection{Setup and preliminaries} \label{S:setup}

We consider the following general form of fractional linear PDE:
\begin{equation}
\label{PDE:w}
q_t+w\left(-i\frac{\partial}{\partial x}\right)q=0,\quad\quad x\in(0,\infty),t\in(0,T),
\end{equation}
where $w$ is a finite fractional series of power functions. More explicitly, we write
\begin{equation}
\label{w:defn}
w(k)=\sum_{\alpha}c_{\alpha}k^{\alpha},\quad\quad \alpha,c_{\alpha}\in\mathbb{C},\mathrm{Re}(\alpha)>0,
\end{equation}
where the summation is finite, i.e. the indices $\alpha$ are contained in some finite set of complex numbers in the right half plane. Thus we can write the PDE more explicitly as
\begin{equation}
\label{PDE:sum}
\frac{\partial q}{\partial t}+\sum_{\alpha}c_{\alpha}(-i)^{\alpha}\frac{\partial^{\alpha}q}{\partial x^{\alpha}}=0,\quad\quad x\in(0,\infty),t\in(0,T).
\end{equation}
We shall attempt to solve this PDE with the initial condition $q(x,0)=q_0(x)$, $x\in\mathbb{R}^+$, where the function $q_0:[0,\infty)\rightarrow\mathbb{C}$ is given and has a well-defined Fourier transform $\hat{q}_0$. We will also need boundary conditions, but their number and nature will be determined later.

For the purposes of this paper, all fractional derivatives are defined in the \textbf{Riemann--Liouville} sense with the constant of differintegration being $0$ -- this is a logical lower bound for integration because the spatial domain for the PDE is bounded below by $0$. We are using the half-line $x\in[0,\infty)$ as our spatial domain because the equation on the full line would be relatively easy to solve using a standard Fourier transform method. For the problem on the half-line, an analytic solution of the PDE requires more advanced transform methods.

We shall be using fractional power functions in the analysis below, and so it will be necessary to define domains and branches for these functions. For our purposes, all fractional power functions are defined using the principal branch with branch cut along the negative real axis, i.e.
\begin{equation}
\label{branch:defn}
k^{\alpha}=r^{\alpha}e^{i\theta\alpha},\quad\quad k=re^{i\theta},r\in\mathbb{R}^+,\theta\in(-\pi,\pi).
\end{equation}

For Fourier transforms with respect to $x$, we will use $\hat{\color{white}q}$ to denote a half-Fourier transform defined on $[0,\infty)$ as follows:
\begin{equation}
\label{xFourier}
\hat{q}(k,t)=\int_0^{\infty}e^{-ikx}q(x,t)\,\mathrm{d}x.
\end{equation}
We do not specify here a particular transform with respect to $t$, because the particular transform we shall use depends on the approach taken, and it will be more complicated than the Fourier transform \eqref{xFourier}.

\subsection{Finding the global relation using a divergence form} \label{S:GRdiv}

Here we work through the first two steps of the method as laid out in \S\ref{intro:UTM}, namely writing the PDE \eqref{PDE:w} in divergence form and deriving a global relation. The biggest challenge here, as we shall see, is to define the function $Q(x,t,k)$ in an appropriate way so that the rest of the argument works.

The fractional PDE \eqref{PDE:w} can still be written in divergence form as \eqref{ex:divergence}, provided that the function $Q(x,t,k)$ satisfies the following condition:
\begin{equation}
\label{Q:condn}
\big(\prescript{}{0}D_x-ik\big)Q=\big(w(k)-w(-i\prescript{}{0}D_x)\big)q.
\end{equation}
However, now that $w$ is no longer necessarily a polynomial, the simple expression \eqref{ex:Qdef} for $Q$ no longer applies. Already the presence of fractional derivatives makes the problem harder than in the classical case. How can we find an explicit form for $Q$ in this case?

Intuitively, we can still consider the function $\frac{w(k)-w(l)}{k-l}$ with the idea of setting $l=-i\prescript{}{0}D_x$ at some later stage. One idea would be to expand $(k-l)^{-1}$ as a power series, namely to use the following:
\[
(k-l)^{-1}=
\begin{cases}
-l^{-1}(1+kl^{-1}+k^2l^{-2}+k^3l^{-3}+\dots)=\sum\limits_{j=0}^{\infty}-k^jl^{-j-1},&\quad|k|<1; \\ \\
k^{-1}(1+k^{-1}l+k^{-2}l^2+k^{-3}l^3+\dots)=\sum\limits_{j=0}^{\infty}k^{-j-1}l^{j},&\quad|k|>1.
\end{cases}
\]
Multiplying this series by $w(k)-w(l)$ would yield an expression for $\frac{w(k)-w(l)}{k-l}$ as an infinite series of terms of the form $k^{\alpha}l^{\beta}$, namely:
\[
i(k-l)^{-1}\left(w(k)-w(l)\right)=
\begin{cases}
i\sum\limits_{j=0}^{\infty}\sum\limits_{\alpha}c_{\alpha}\left[-k^{\alpha+j}l^{-j-1}+k^{j}l^{\alpha-j-1}\right],&\quad|k|<1; \\ \\
i\sum\limits_{j=0}^{\infty}\sum\limits_{\alpha}c_{\alpha}\left[k^{\alpha-j-1}l^{j}-k^{-j-1}l^{\alpha+j}\right],&\quad|k|>1.
\end{cases}
\]
However, fractional differential operators do not have a semigroup property, by Lemma \ref{ComposeDerivatives}: after setting $l=-i\prescript{}{0}D_x$, the product of $l^a$ and $l^b$ will not necessarily be $l^{a+b}$. So the above manipulation of terms is actually not valid if $l=-i\prescript{}{0}D_x$ is assumed \textit{a priori}.

Fortunately, we do not need $Q$ to be precisely the expression $i\left(\frac{w(k)-w(l)}{k-l}\right)$ with $l$ replaced by $-i\prescript{}{0}D_x$. Any function $Q$ that satisfies the condition \eqref{Q:condn} will automatically give the divergence form \eqref{ex:divergence} as an equivalent formulation of the PDE \eqref{PDE:w}. So we only need to find a function $Q$ satisfying \eqref{Q:condn}, and considering $\frac{w(k)-w(l)}{k-l}$ provides a motivation for where to look for such a function. With this in mind, let us try defining $Q$ by the above series with $l=-i\prescript{}{0}D_x$, ignoring whether our manipulations of $l$ would actually be valid for this differential operator. The resulting expression for $Q$ is:
\begin{equation}
\label{Q:defn}
Q(x,t,k)=
\begin{cases}
i\sum\limits_{j=0}^{\infty}\sum\limits_{\alpha}c_{\alpha}\left[k^{j}\left(-i\prescript{}{0}D_x\right)^{\alpha-j-1}(q)-k^{\alpha+j}\left(-i\prescript{}{0}D_x\right)^{-j-1}(q)\right],&\quad|k|<1; \\ \\
i\sum\limits_{j=0}^{\infty}\sum\limits_{\alpha}c_{\alpha}\left[k^{\alpha-j-1}\left(-i\prescript{}{0}D_x\right)^{j}(q)-k^{-j-1}\left(-i\prescript{}{0}D_x\right)^{\alpha+j}(q)\right] + A(k)e^{ikx},&\quad|k|>1;
\end{cases}
\end{equation}
where $A(k)$ is chosen so that $Q$ is continuous across $|k|=1$. Such a function $A(k)$ exists because $A(k)e^{ikx}$ is the general solution of the ODE $\big(\prescript{}{0}D_x-ik\big)Q=0$ in $x$.

It is now straightforward to verify, using the fact that $\prescript{}{0}D_x\circ\prescript{}{0}D_x^{\alpha}=\prescript{}{0}D_x^{\alpha+1}$ for all $\alpha\in\mathbb{C}$ by Lemma \ref{ComposeDerivatives}, that both of the above series expressions for $Q(x,t,k)$ do satisfy \eqref{Q:condn} as required. So \eqref{Q:defn} provides a possible choice for the function $Q$.

The problem now is that an infinite series expression for $Q$ is difficult to deal with. It looks as though integrating the divergence form \eqref{ex:divergence} and using the expressions \eqref{Q:defn} for $Q$ will lead to a global relation in the form of an infinite series, with infinitely many boundary terms involved, and therefore requiring infinitely many boundary conditions specified in order to have a unique solution. However, it is already known (see e.g. \cite{kaikina2}) that only finitely many boundary conditions need to be specified in order to get a unique solution for a PDE of this form.

But in fact it turns out that almost all the terms in the infinite series for $Q$ cancel out when substituted into the global relation!

The next step in the method as laid out in \S\ref{intro:UTM} is to apply integrals with respect to both $x$ and $t$ to the divergence form \eqref{ex:divergence} of the PDE. This yields the following relation between integral transforms:
\begin{equation}
\label{GR:first}
e^{w(k)t}\hat{q}(k,t)-\hat{q}_0(k)=-\int_0^te^{w(k)\tau}Q(0,\tau,k)d\tau,\quad\quad\mathrm{Im}(k)<0.
\end{equation}
Here we have assumed sufficient decay conditions that the upper limit term for $x$ vanishes, i.e.
\begin{equation}
\label{decay}
\lim_{x\rightarrow\infty}\int_0^te^{-ikx+w(k)\tau}Q(x,\tau,k)d\tau=0.
\end{equation}

Importantly, the only time $Q$ appears in \eqref{GR:first} is when $x=0$: we do not need to deal with the full complexity of the function $Q(x,t,k)$, but only with the special case $Q(0,\tau,k)$. And by definition of the Riemann--Liouville fractional integral, the function $\prescript{}{0}D_x^{\nu}f(x)\big|_{x=0}$ for any given $f$ is always identically zero when this is a fractional integral, i.e. when $\mathrm{Re}(\nu)<0$. So for the purposes of the global relation \eqref{GR:first}, we can ignore all terms in the infinite series of \eqref{Q:defn} in which $\prescript{}{0}D_x$ appears to a negative power.

Thus, we only consider the $|k|<1$ part of \eqref{Q:defn}, since negative powers of $\prescript{}{0}D_x$ appear in only finitely many terms of the $|k|>1$ part and in all but finitely many terms of the $|k|<1$ part. Substituting $x=0$ into the $|k|<1$ series from \eqref{Q:defn}, we find:
\begin{align*}
Q(0,\tau,k)&=i\sum_{j=0}^{\infty}\sum_{\alpha}c_{\alpha}k^{j}\left(-i\prescript{}{0}D_x\right)^{\alpha-j-1}q(0,\tau) \\
&=-\sum_{\alpha}c_{\alpha}(-i)^{\alpha}\sum_{j=0}^{\lfloor\mathrm{Re}(\alpha)\rfloor-1}(ik)^{j}\prescript{}{0}D_x^{\alpha-j-1}q(0,\tau),\quad\quad|k|<1.
\end{align*}
Thus the identity \eqref{GR:first} yields the following global relation:
\begin{equation}
\label{GR:final}
e^{w(k)t}\hat{q}(k,t)-\hat{q}_0(k)=\sum_{\alpha}c_{\alpha}(-i)^{\alpha}\sum_{j=0}^{\lfloor\mathrm{Re}(\alpha)\rfloor-1}(ik)^{j}\int_0^te^{w(k)\tau}\prescript{}{0}D_x^{\alpha-j-1}q(0,\tau)d\tau,
\end{equation}
valid for $\mathrm{Im}(k)<0$ and $|k|<1$. Now that the infinite series over $j$ has become a finite one, we no longer need $|k|<1$ for convergence. So by analytic continuation, \eqref{GR:final} is valid for all $k$ in the lower half plane $\mathrm{Im}(k)<0$. Thus, \eqref{GR:final} provides a finite closed-form global relation as desired.

\subsection{Finding the global relation using double transforms} \label{S:GRfour}

In this section, we consider another way of deriving the global relation. This does not follow the approach which was indicated in \S\ref{intro:UTM} for non-fractional PDEs, but it is similar to a known alternative methodology \cite{arciga4}, and it yields a global relation equivalent to the one found in \S\ref{S:GRdiv}.

Before proceeding to analyse the PDE, we first obtain an identity which we shall need to use in this section. Applying the fractional integration by parts rule \eqref{IbP:fractional} to the functions $f(x)=e^{-ikx}$ and $g(x)=q(x,t)$, with upper and lower limits $a=0$ and $b\rightarrow\infty$, and using the formula \eqref{RL:exp} for fractional differintegrals of exponential functions, we find
\begin{multline*}
\int_0^{\infty}e^{-ikx}\cdot \prescript{}{0}D_{x}^{\alpha}q(x,t)\,\mathrm{d}x \\ =\int_0^{\infty}q(x,t)\cdot(ik)^{\alpha}e^{-ikx}\,\mathrm{d}x-\sum_{j=0}^{n-1}\left[(-1)^{n+j}\prescript{}{0}D_{x}^{\alpha-n+j}q(x,t)\cdot(-ik)^{n-j-1}e^{-ikx}\right]_0^{\infty},
\end{multline*}
or in other words
\begin{equation}
\label{IbP:q}
\widehat{\frac{\partial^{\alpha}q}{\partial x^{\alpha}}}(k,t)=(ik)^{\alpha}\hat{q}(k,t)+\sum_{j=0}^{n-1}(ik)^{n-j-1}\left[e^{-ikx}\prescript{}{0}D_{x}^{\alpha-n+j}q(x,t)\right]_0^{\infty}.
\end{equation}

We also use $\tilde{•}$ to denote the Laplace transform with respect to $t$:
\begin{equation}
\label{tLaplace}
\tilde{q}(x,s)=\int_0^{\infty}e^{-st}q(x,t)\,\mathrm{d}t.
\end{equation}
For the purposes of the definition \eqref{tLaplace}, we extend the function $q$ beyond the interval $[0,T]$ by making it identically zero for large $t$. This will have no effect on the final result of this section, because the only appearance of Laplace transforms will be to be applied and then almost immediately inverted again.

Armed with the integration by parts identity \eqref{IbP:q}, we proceed to apply a half-Fourier transform with respect to $x$ to the PDE \eqref{PDE:sum}:
\begin{align*}
\eqref{PDE:sum}&\Rightarrow\frac{\partial\hat{q}}{\partial t}+\sum_{\alpha}c_{\alpha}(-i)^{\alpha}\widehat{\frac{\partial^{\alpha}q}{\partial x^{\alpha}}}=0 \\
&\Rightarrow\frac{\partial\hat{q}}{\partial t}+\sum_{\alpha}c_{\alpha}(-i)^{\alpha}\left[(ik)^{\alpha}\hat{q}(k,t)+\sum_{j=0}^{n-1}(ik)^{n-j-1}\left[e^{-ikx}\prescript{}{0}D_{x}^{\alpha-n+j}q(x,t)\right]_0^{\infty}\right]=0 \\
&\Rightarrow\frac{\partial\hat{q}}{\partial t}+w(k)\hat{q}+\sum_{\alpha}c_{\alpha}(-i)^{\alpha}\sum_{j=0}^{n-1}(ik)^{n-j-1}\left[e^{-ikx}\prescript{}{0}D_{x}^{\alpha-n+j}q(x,t)\right]_0^{\infty}=0.
\end{align*}
Note that the notation $n=\lceil\mathrm{Re}(\alpha)\rceil$ is propagated here from Lemma \ref{IbPlemma}.

Next, we apply a Laplace transform with respect to $t$. This results in an expression for the double-transformed function $\tilde{\hat{q}}(k,s)$: \[s\tilde{\hat{q}}(k,s)-\hat{q}(k,0)+w(k)\tilde{\hat{q}}(k,s)+\sum_{\alpha}c_{\alpha}(-i)^{\alpha}\sum_{j=0}^{n-1}(ik)^{n-j-1}\left[e^{-ikx}\prescript{}{0}D_{x}^{\alpha-n+j}\tilde{q}(x,s)\right]_0^{\infty}=0,\] which rearranges to
\begin{equation*}
\tilde{\hat{q}}(k,s)=\frac{\hat{q}(k,0)+\sum_{\alpha}c_{\alpha}(-i)^{\alpha}\sum_{j=0}^{n-1}(ik)^{n-j-1}\left[e^{-ikx}\prescript{}{0}D_{x}^{\alpha-n+j}\tilde{q}(x,s)\right]_0^{\infty}}{s+w(k)}.
\end{equation*}
Now we have an explicit formula, if not for $q$ itself, at least for some transform of $q$, in terms of its initial and boundary values. We apply an inverse Laplace transform with respect to $t$, recalling both the convolution theorem and the fact that the transform of $e^{-w(k)t}$ is $\frac{1}{s+w(k)}$:
\begin{align*}
\hat{q}(k,t)&=\hat{q}(k,0)e^{-w(k)t}-\left[\sum_{\alpha}c_{\alpha}(-i)^{\alpha}\sum_{j=0}^{n-1}(ik)^{n-j-1}\left[e^{-ikx}\prescript{}{0}D_{x}^{\alpha-n+j}q(x,t)\right]_0^{\infty}\right]\ast\left[e^{-w(k)t}\right] \\
&=\hat{q}(k,0)e^{-w(k)t}-\int_0^t\sum_{\alpha}c_{\alpha}(-i)^{\alpha}\sum_{j=0}^{n-1}(ik)^{n-j-1}\left[e^{-ikx}\prescript{}{0}D_{x}^{\alpha-n+j}q(x,\tau)\right]_0^{\infty}e^{-w(k)(t-\tau)}\,\mathrm{d}\tau.
\end{align*}
Thus the global relation is
\begin{equation}
\label{GR}
e^{w(k)t}\hat{q}(k,t)=\hat{q}(k,0)-\sum_{\alpha}c_{\alpha}(-i)^{\alpha}\sum_{j=0}^{n-1}(ik)^{n-j-1}\left[e^{-ikx}\int_0^te^{w(k)\tau}\prescript{}{0}D_{x}^{\alpha-n+j}q(x,\tau)\,\mathrm{d}\tau\right]_0^{\infty}.
\end{equation}
It is valid for $\mathrm{Im}(k)<0$, this condition being required for the half-Fourier transforms with respect to $x$ to be well-defined.

For ease of notation, we define
\begin{equation}
\label{g:defn}
g(k,t):=\sum_{\alpha}c_{\alpha}(-i)^{\alpha}\sum_{j=0}^{n-1}(ik)^{n-j-1}\left[e^{-ikx}\int_0^te^{w(k)\tau}\prescript{}{0}D_{x}^{\alpha-n+j}q(x,\tau)\,\mathrm{d}\tau\right]_0^{\infty},
\end{equation}
so that the global relation is
\begin{equation}
\label{GR:g}
e^{w(k)t}\hat{q}(k,t)=\hat{q}(k,0)-g(k,t),\quad\quad\mathrm{Im}(k)<0.
\end{equation}

We note that, as expected, the global relation \eqref{GR} is identical, under the assumption \eqref{decay} on the decay of $q$ at infinity, to the previously obtained global relation \eqref{GR:final}. The discrepancy in the number of terms in the series (namely, $\lfloor\mathrm{Re}(\alpha)\rfloor$ in \eqref{GR:final} versus $\lceil\mathrm{Re}(\alpha)\rceil$ in \eqref{GR}) is resolved by using the same argument as in \S\ref{S:GRdiv} to point out that the fractional integral $\prescript{}{0}D_{x}^{\alpha-\lceil\mathrm{Re}(\alpha)\rceil}q(0,t)$ is identically zero for any $\alpha$ with non-integer real part.

In summary, we have obtained exactly the same identity twice using two different approaches.

\subsection{Deducing the solution} \label{S:laststep}

We start from the global relation \eqref{GR:g} and apply an inverse half-Fourier transform with respect to $x$. This yields the following explicit expression for $q$:
\begin{equation}
\label{q:first}
q(x,t)=\frac{1}{2\pi}\int_{-\infty}^{\infty}e^{ikx-w(k)t}\left(\hat{q}(k,0)-g(k,t)\right)\,\mathrm{d}k.
\end{equation}
Note that the function $e^{ikx-w(k)t}g(k,t)$ is analytic in $k$ everywhere except along the branch cut for $w(k)$; so by our definition \eqref{branch:defn}, it is analytic on the domain $\mathbb{C}\backslash(-\infty,0]$ for $k$. Furthermore it has exponential decay (tends to zero) as $|k|\rightarrow\infty$ with $\mathrm{Im}(k)>0,\mathrm{Re}(w(k))>0$. So by Cauchy's theorem, we can deform the contour of integration for the second half of the integral in \eqref{q:first} through any region with $\mathrm{Im}(k)>0$ and $\mathrm{Re}(w(k))>0$. Thus, we obtain the following improved explicit formula for $q$:
\begin{equation}
\label{q:second}
q(x,t)=\frac{1}{2\pi}\int_{-\infty}^{\infty}e^{ikx-w(k)t}\hat{q}(k,0)\,\mathrm{d}k-\frac{1}{2\pi}\int_{\partial D^+}e^{ikx-w(k)t}g(k,t)\,\mathrm{d}k,
\end{equation}
where the domain $D^+$ is defined as
\begin{equation}
\label{D+:defn}
D^+=\{k\in\mathbb{C}:\mathrm{Im}(k)>0,\mathrm{Re}(w(k))<0\}
\end{equation}
and we assume that $\mathrm{Re}(w(k))\geq0$ for all $k\in\mathbb{R}$.

The equation \eqref{q:second} gives us an expression for $q$ in terms of a single initial condition, namely $q(x,0)$, and $n$ boundary values, namely $D_x^{\alpha-n+j}q(0,t)$ for $0\leq j<n-1$, provided we have sufficient decay conditions on $q$ as $x\rightarrow\infty$. But this is an overdetermined problem, since we cannot prescribe these $n$ boundary values as boundary conditions \cite{kaikina1,kaikina2}.

Fortunately, it is possible to eliminate some of the boundary values by using the global relation \eqref{GR} with $k$ replaced by $\nu(k)$ for some function $\nu$. This function is required to satisfy two properties:
\begin{itemize}
\item It must preserve $w$, i.e. $w(\nu(k))=w(k)$. This is so that all the terms involving $w(k)$ in \eqref{GR} are not altered by the substitution, while those directly involving $k$ may change.
\item It must map the boundary $\partial D^+$ into some region of the lower half $k$-plane. This is so that the global relation is valid at $\nu(k)$ when $k\in\partial D^+$.
\end{itemize}
Such functions $\nu$ are hard to find explicitly for the most general function $w$. But given a specific $w$, it is often possible to find the $\nu$ required. For example, in the simple case of $w(k)=k^{\alpha}$, we can take $\nu(k)=e^{2\pi im/\alpha}$ for any integer value of $m$; this problem is considered in detail in \S\ref{S+:example} below.

Using the functions $\nu$, we find new equations in a similar form to \eqref{GR} which are valid for $k\in\partial D^+$ and therefore can be used in \eqref{q:second}. By making multiple such substitutions, it is possible to eliminate the unknown boundary values.

\subsection{Summary and verification} \label{S:verify}
Following the above steps, we obtain the following result.

\begin{thm}
Given a PDE of the form \eqref{PDE:w} valid on the region $0<x<\infty,0<t<T$, where $w$ is a finite series of power functions defined by \eqref{w:defn} and satisfying $\mathrm{Re}(w(k))\geq0$ for all $k\in\mathbb{R}$, and where fractional derivatives are defined in the Riemann--Liouville sense \eqref{RLdef:deriv} with $a=0$, the unified transform method can be used to construct an explicit solution $q(x,t)$ in terms of the initial condition $q(x,0)$ and some boundary conditions $\frac{\partial^{\alpha-r}q}{\partial x^{\alpha-r}}(0,t)$.
\end{thm}

In order to verify this result, we substitute the formula \eqref{q:second} into the original PDE \eqref{PDE:w} to check that this $q$ does indeed satisfy the equation. (The final formula would be in a more complicated form than \eqref{q:second}, but from that formula it is easy to return to the expression \eqref{q:second} just by reversing the substitutions made to get there -- all modifications between the two are only a matter of rewriting boundary values in terms of each other.)

Starting from the formula \eqref{q:second} for $q$, we find
\begin{multline*}
\frac{\partial q}{\partial t}=-\frac{1}{2\pi}\int_{-\infty}^{\infty}w(k)e^{ikx-w(k)t}\hat{q}(k,0)\,\mathrm{d}k+\frac{1}{2\pi}\int_{\partial D^+}w(k)e^{ikx-w(k)t}g(k,t)\,\mathrm{d}k \\ -\frac{1}{2\pi}\int_{\partial D^+}e^{ikx-w(k)t}g_t(k,t)\,\mathrm{d}k,
\end{multline*}
and, using \eqref{IbP:q},
\begin{equation*}
w\left(-i\frac{\partial}{\partial x}\right)q=\frac{1}{2\pi}\int_{-\infty}^{\infty}e^{ikx}\left[w(k)\hat{q}(k,t)+\sum_{\alpha}c_{\alpha}(-i)^{\alpha}\sum_{j=0}^{n-1}(ik)^{n-j-1}\left[e^{-ikx}\prescript{}{0}D_{x}^{\alpha-n+j}q(x,t)\right]_0^{\infty}\right]\,\mathrm{d}k.
\end{equation*}
So the left-hand side of the original PDE evaluates as follows:
\begin{multline*}
\frac{\partial q}{\partial t}+w\left(-i\frac{\partial}{\partial x}\right)q \\ =-\frac{1}{2\pi}\int_{-\infty}^{\infty}w(k)e^{ikx-w(k)t}\hat{q}(k,0)\,\mathrm{d}k+\frac{1}{2\pi}\int_{\partial D^+}w(k)e^{ikx-w(k)t}g(k,t)\,\mathrm{d}k \\ \hspace{2cm}-\frac{1}{2\pi}\int_{\partial D^+}e^{ikx-w(k)t}\sum_{\alpha}c_{\alpha}(-i)^{\alpha}\sum_{j=0}^{n-1}(ik)^{n-j-1}\left[e^{-ikx+w(k)t}\prescript{}{0}D_{x}^{\alpha-n+j}q(x,t)\right]_0^{\infty}\,\mathrm{d}k \\ +\frac{1}{2\pi}\int_{-\infty}^{\infty}w(k)e^{ikx}\left[e^{-w(k)t}\hat{q}(k,0)-e^{-w(k)t}g(k,t)\right]\mathrm{d}k \\ \hspace{-3cm}+\frac{1}{2\pi}\int_{-\infty}^{\infty}e^{ikx}\sum_{\alpha}c_{\alpha}(-i)^{\alpha}\sum_{j=0}^{n-1}(ik)^{n-j-1}\left[e^{-ikx}\prescript{}{0}D_{x}^{\alpha-n+j}q(x,t)\right]_0^{\infty}\,\mathrm{d}k,
\end{multline*}
which is zero since Cauchy's theorem enables us to equate $\int_{\partial D^+}$ with $\int_{-\infty}^{\infty}$ when required. We also used the global relation \eqref{GR} as a substitution for $\hat{q}(k,t)$ in the derivation of the above formula.

Thus, we have proved that the solution constructed above actually is a solution, subject to a straightforward verification of the initial and boundary conditions.

\section{Applications and extensions} \label{S+:example}

\subsection{A worked example}

As a basic but important example of the method outlined in the previous section, let us consider the case where $w$ is a single power function, say \[w(k)=-A(ik)^{\alpha},\] where the constants $A$ and $\alpha$ are fixed. For simplicity we now assume both of these constants to be positive real.

In other words, we shall attempt to solve the following PDE:
\begin{equation}
\label{PDE:example}
\frac{\partial q}{\partial t}=A\frac{\partial^{\alpha}q}{\partial x^{\alpha}},\quad\quad x\in(0,\infty),t\in(0,T).
\end{equation}
As always, we impose a single initial condition, \[q(x,0)=q_0(x),\quad\quad x\in(0,\infty),\] and a certain number of boundary conditions to be determined later.

In this case, the global relation \eqref{GR:final} is:
\begin{equation}
\label{GR:example}
e^{-A(ik)^{\alpha}t}\hat{q}(k,t)=\hat{q}_0(k)-A\sum_{j=0}^{\lfloor\alpha\rfloor-1}(ik)^{j}\int_0^te^{-A(ik)^{\alpha}\tau}\prescript{}{0}D_x^{\alpha-j-1}q(0,\tau)d\tau,\quad\quad\mathrm{Im}(k)<0.
\end{equation}

In order to find the region $D^+$ defined by \eqref{D+:defn}, we note that the following conditions are equivalent:
\begin{align*}
&\mathrm{Re}(w(k))<0; \\
&\mathrm{Re}\left((ik)^{\alpha}\right)>0; \\
2n\pi-\frac{\pi}{2}<&\arg\left((ik)^{\alpha}\right)<2n\pi-\frac{\pi}{2}\quad\text{ for some }n\in\mathbb{Z}; \\
\frac{\pi}{\alpha}\left(2n-\frac{1}{2}\right)<&\arg(ik)<\frac{\pi}{\alpha}\left(2n+\frac{1}{2}\right)\quad\text{ for some }n\in\mathbb{Z}.
\end{align*}
We saw in \S\ref{S:laststep} that all real $k$ are required \textit{not} to satisfy this condition, in order that a meaningful deformation of contours can be applied. Thus we require that \[\left|\frac{2n\pi}{\alpha}\pm\frac{\pi}{2}\right|\geq\frac{\pi}{2\alpha}\] for all integers $n$, i.e. that $|4n\pm\alpha|\geq1$ for all integers $n$. In other words, $\alpha$ must lie in one of the intervals $[1,3]$, $[5,7]$, $[9,11]$, etc.

The domain $D^+$ can be described as follows, working from \eqref{D+:defn}:
\begin{align*}
D^+&=\{k\in\mathbb{C}:\mathrm{Im}(k)>0,\mathrm{Re}\left((ik)^{\alpha}\right)>0\} \\
&=\Big\{k\in\mathbb{C}:0<\arg(k)<\tfrac{\pi}{2},2n\pi-\tfrac{\pi}{2}<\alpha\left[\arg(k)+\tfrac{\pi}{2}\right]<2n\pi+\tfrac{\pi}{2},n\in\mathbb{Z}\Big\} \\
&\hspace{2cm}\cup\Big\{k\in\mathbb{C}:\tfrac{\pi}{2}<\arg(k)<\pi,2n\pi-\tfrac{\pi}{2}<\alpha\left[\arg(k)-\tfrac{3\pi}{2}\right]<2n\pi+\tfrac{\pi}{2},n\in\mathbb{Z}\Big\} \\
&=\Big\{k\in\mathbb{C}:\alpha\left[\arg(k)+\tfrac{\pi}{2}\right]\in\left(\tfrac{\pi\alpha}{2},\pi\alpha\right)\cap\left(2n\pi-\tfrac{\pi}{2},2n\pi+\tfrac{\pi}{2}\right),n\in\mathbb{Z}\Big\} \\
&\hspace{2cm}\cup\Big\{k\in\mathbb{C}:\alpha\left[\arg(k)-\tfrac{3\pi}{2}\right]\in\left(-\pi\alpha,-\tfrac{\pi\alpha}{2}\right)\cap\left(2n\pi-\tfrac{\pi}{2},2n\pi+\tfrac{\pi}{2}\right),n\in\mathbb{Z}\Big\}.
\end{align*}
In particular, if $1<\alpha\leq\frac{3}{2}$, then the intervals on the right do not intersect and so $D^+$ is empty. For now, let us assume $\boldsymbol{\frac{3}{2}<\alpha<\frac{5}{2}}$, so that we have:
\begin{align}
\nonumber D^+&=\Big\{k\in\mathbb{C}:\alpha\left[\arg(k)+\tfrac{\pi}{2}\right]\in\left(\tfrac{\pi\alpha}{2},\pi\alpha\right)\cap\left(\tfrac{3\pi}{2},\tfrac{5\pi}{2}\right)\Big\} \\
\nonumber &\hspace{2cm}\cup\Big\{k\in\mathbb{C}:\alpha\left[\arg(k)-\tfrac{3\pi}{2}\right]\in\left(-\pi\alpha,-\tfrac{\pi\alpha}{2}\right)\cap\left(-\tfrac{5\pi}{2},-\tfrac{3\pi}{2}\right)\Big\} \\
\label{D+:example} &=\Big\{k\in\mathbb{C}:\tfrac{3\pi}{2\alpha}-\tfrac{\pi}{2}<\arg(k)<-\tfrac{3\pi}{2\alpha}+\tfrac{3\pi}{2}\Big\}.
\end{align}

The integral expression \eqref{q:second} for $q$ now becomes:
\begin{equation}
\label{q:example}
q(x,t)=\frac{1}{2\pi}\int_{-\infty}^{\infty}e^{ikx+A(ik)^{\alpha}t}\hat{q_0}(k)\,\mathrm{d}k-\frac{1}{2\pi}\int_{\Gamma}e^{ikx+A(ik)^{\alpha}t}g(k,t)\,\mathrm{d}k,
\end{equation}
with notation defined as follows. The contour $\Gamma$ runs along the boundary of $D^+$, forming a V-shape in the upper half plane:
\begin{equation}
\label{Gamma:defn}
\Gamma=\left\{re^{(3\alpha-3)\pi i/2\alpha} : \infty>r>0\right\}\cup\left\{re^{(3-\alpha)\pi i/2\alpha} : 0<r<\infty\right\}.
\end{equation}
And the function $g(k,t)$ is defined by minus the main term on the right-hand side of \eqref{GR:example}, this being the notation required for the global relation to be expressible in the form \eqref{GR:g}. In other words, $g$ is given by the following expression:
\begin{equation}
\label{g:example}
g(k,t)=A\sum_{j=0}^{\lfloor\alpha\rfloor-1}(ik)^{j}\int_0^te^{-A(ik)^{\alpha}\tau}\prescript{}{0}D_x^{\alpha-j-1}(q)d\tau.
\end{equation}

Now we need to find functions $\nu$ which preserve the power function $w$, i.e. such that \[(\nu(k))^{\alpha}=k^{\alpha}.\] This is easy to solve for $\nu$; the function \[\nu(k)=e^{2n\pi i/\alpha}k\] will work for any integer $n$ such that $\arg(k)$ and $\frac{2n\pi}{\alpha}+\arg(k)$ are both in the domain $(-\pi,\pi)$ required by the power function with branch cut along the negative real axis. We also require $\frac{2n\pi}{\alpha}+\arg(k)$ to be in the domain $(-\pi,0)$ when $k$ is on the contour $\Gamma$, in order that our substitution into the global relation will be valid.

Clearly any positive $n$ would not transform the contour $\Gamma$ into the lower half plane. Furthermore, any $n\leq-2$ would transform $\Gamma$ outside of the domain $(-\pi,\pi)$ for arguments. The case $n=0$ only gives us the identity map. So the only non-trivial possibility for $\nu$ is with $n=-1$, namely:
\begin{equation}
\label{nu:example}
\nu(k)=e^{-2\pi i/\alpha}k.
\end{equation}
This function $\nu$ transforms $\Gamma$ into the lower half plane if and only if $\alpha>\frac{7}{3}$. Let us now assume $\boldsymbol{2<\alpha<\frac{7}{3}}$, so that we have both a valid map $\nu$ and a fixed value of $\lfloor\alpha\rfloor$. Now, substituing $k$ for $\nu(k)$ into the global relation \eqref{GR:example}, we find:
\begin{multline}
\label{transformedGR}
e^{-A(ik)^{\alpha}t}\hat{q}\left(e^{-2\pi i/\alpha}k,t\right) \\ =\hat{q}_0\left(e^{-2\pi i/\alpha}k\right)-A\sum_{j=0}^{1}\left(ie^{-2\pi i/\alpha}k\right)^{j}\int_0^te^{-A(ik)^{\alpha}\tau}\prescript{}{0}D_x^{\alpha-j-1}q(0,\tau)d\tau,\quad\quad k\in\Gamma.
\end{multline}
Substituting \eqref{transformedGR} into \eqref{q:example} enables us to eliminate one of the $\lfloor\alpha\rfloor=2$ boundary conditions on the right-hand side, leaving only one boundary condition that needs to be specified in the initial setup of the problem.

Thus, in the case $2<\alpha<\frac{7}{3}$, the unified transform method can be used to solve the fractional PDE \eqref{PDE:example}, with the initial condition $q(x,0)=q_0(x)$ and exactly one of the two boundary terms $\prescript{}{0}D_x^{\alpha-1}q(0,t)$, $\prescript{}{0}D_x^{\alpha-2}q(0,t)$ specified.

Of course, this range of values of $\alpha$ is not the only one in which the problem \eqref{PDE:example} can be solved. We chose our restrictions on $\alpha$ merely for convenience. It would be just as easy to solve the PDE in the case $\frac{5}{2}<\alpha<3$, or other higher ranges of $\alpha$. What we have presented here is the solution of a model problem, in order to demonstrate the methodology. Other example problems would work out similarly, but might become more complicated according to the value of $\alpha$.

\subsection{Potential extensions}

There are many ways in which the method laid out in \S\ref{S+:main} could be generalised beyond even the general equation \eqref{PDE:w}.

In the example discussed above, we assumed that the index $\alpha$ was real. Even the simple PDE \eqref{PDE:example} becomes more interesting to solve when $\alpha$ is complex. The boundary of the domain $D^+$ would no longer consist of rays from the origin in the complex plane, but rather of infinite logarithmic spirals (due to considering the argument of $k^{\alpha}$ with $\alpha$ complex), and the contour of integration would become correspondingly more complicated.

For integer-order PDEs, the unified transform method has been applied to many families of equations more advanced than \eqref{ex:PDE} on more complicated domains than the half-line $[0,\infty)$ -- for example, finite intervals in the real line, convex polygons in a plane, and beyond \cite{fokas3}. Fractional analogues of these problems could be considered and potentially solved by modifying the unified transform method for the new scenario.

In this paper we have considered only fractional differential equations of Riemann--Liouville type. But many real-world processes can be better modelled using other definitions of fractional calculus: for example, the Caputo definition \eqref{CAPdef} is better suited to many initial value problems, and newer definitions such as Caputo--Fabrizio and Atangana--Baleanu have been used to model various types of nonlocal dynamics \cite{atangana2,atangana,hristov}. Solving fractional PDEs in these alternative fractional models could be an important result, and it may be possible to do so using the unified transform method. The Riemann--Liouville model has the unique advantage of interacting with Fourier and Laplace transforms in the way one would expect for a fractional derivative, namely with fractional derivatives becoming fractional power function multipliers in the transformed space, but other models have similar properties with the power function replaced by more complicated functions, and PDEs in these models are still amenable to transform approaches.

\section{Conclusions}

In this paper, we have considered the unified transform method as it applies to linear evolution PDEs of fractional order. Various problems arose due to the fractionalisation of the problem, but we demonstrated how each of these problems could be overcome by the introduction of new ideas and methods. We described the method as it applies to a general form of PDE, and then demonstrated its applicability by using it to solve a specific model problem. We also considered several directions in which our work here could be extended in the future.

\section*{Acknowledgements}
The first and third authors acknowledge the support of the Engineering and Physical Sciences Research Council, UK: the first author via a research studentship (no grant number) and the third author via a senior fellowship (grant number RG81609).

\end{document}